\documentclass[12pt]{article}
\usepackage{amsfonts}
\usepackage{mathrsfs}
\usepackage{amsmath,amsfonts,amssymb,rotating,graphics,ifpdf}

\def\qed{\nopagebreak\hfill{\rule{4pt}{7pt}}}
\def\proof{\noindent {\it{Proof.} \hskip 2pt}}

\newtheorem{theorem}{Theorem}[section]

\newtheorem{corollary}[theorem]{Corollary}
\newtheorem{conjecture}[theorem]{Conjecture}

\newdimen\Squaresize \Squaresize=11pt
\newdimen\Thickness \Thickness=0.7pt
\def\Square#1{\hbox{\vrule width \Thickness
   \vbox to \Squaresize{\hrule height \Thickness\vss
    \hbox to \Squaresize{\hss#1\hss}
   \vss\hrule height\Thickness}
\unskip\vrule width \Thickness} \kern-\Thickness}

\def\Vsquare#1{\vbox{\Square{$#1$}}\kern-\Thickness}

\def\moins{\raise 1pt\hbox{{$\scriptstyle -$}}}

\usepackage{color}

\begin{document}

\begin{center}
{\large \bf  On  Balanced Colorings of the $n$-Cube}
\end{center}

\begin{center}
William Y.C. Chen$^{1}$ and Larry X.W. Wang$^{3}$\\[6pt]
Center for Combinatorics, LPMC-TJKLC\\
Nankai University, Tianjin 300071, P. R. China\\
Email: $^{1}${\tt chen@nankai.edu.cn},  $^{3}${\tt
wxw@cfc.nankai.edu.cn}
\end{center}

\vspace{0.3cm} \noindent{\bf Abstract.}
 A $2$-coloring of the
$n$-cube in the  $n$-dimensional Euclidean space  can be considered
as an assignment of weights of 1 or 0 to the vertices. Such a
colored
 $n$-cube is said to be balanced if its center of mass coincides
with its geometric center. Let $B_{n,2k}$ be the number of
balanced $2$-colorings of the $n$-cube with $2k$ vertices having
weight 1. Palmer, Read and Robinson conjectured that for $n\geq 1$,
the sequence $\{B_{n,2k}\}_{k=0, 1,  \ldots, 2^{n-1}}$ is
symmetric and unimodal. We give a proof of this
conjecture. We also propose a conjecture on the log-concavity of
$B_{n,2k}$ for fixed $k$, and by  probabilistic method we show that
 it holds  when $n$ is sufficiently large.

\noindent {\bf Keywords:} unimodalily, $n$-cube, balanced
coloring

\noindent {\bf AMS Classification:} 05A20, 05D40

\noindent {\bf Suggested Running Title:} Balanced Colorings of the
$n$-Cube

\section{Introduction}

This paper is concerned with a conjecture of Palmer, Read and
Robinson \cite{palmer}  in the $n$-dimensional Euclidean space. A
$2$-coloring of the $n$-cube is considered as an assignment of weights
of 1 or 0 to the vertices. The black vertices are considered as
having weight $1$ whereas  the white vertices are considered as
having weight $0$. We say that a $2$-coloring of the $n$-cube is
balanced if the colored $n$-cube is balanced, namely, the center of
mass is  located at its geometric center.

 Let  $\mathcal{B}_{n,2k}$  denote the set of balanced $2$-colorings
of the $n$-cube with exactly $2k$ black vertices and
$B_{n,2k}=|\mathcal{B}_{n,2k}|$.  Palmer, Read and Robinson proposed
the conjecture that the sequence $\{B_{n,2k}\}_{1\leq k\leq 2^n}$ is
unimodal with the maximum at $k=2^{n-1}$ for any $n\geq 1$. For
example, when $n=4$, the sequence $\{B_{n,2k}\}$ reads
\[
1,8,52,152,222,152,52,8,1.
\]

A sequence $\{a_i\}_{0\leq i\leq m}$ is called unimodal  if there
exists $k$ such that
\[
a_0\leq \cdots \leq a_k\geq  \cdots \geq a_m,
\]
and is called strictly unimodal if
\[
a_0<\cdots < a_k>  \cdots > a_m.
\]
A sequence $\{a_i\}_{0\leq i\leq m}$ of real numbers is said to be
log-concave if
\[
a_i^2\geq a_{i+1}a_{i-1}
\]
 for all $1\leq i\leq m-1$.

Palmer, Read and Robinson \cite{palmer} used P\'{o}lya's theorem to
derive a formula for ${B}_{n,2k}$, which is a sum over integer
partitions of $2k$. However, the unimodality of the sequence
$\{B_{n, 2k}\}$ does not seem to be an easy consequence since the
summation involves negative terms. In Section 2, we will establish a
relation on a refinement of the numbers $\mathcal{B}_{n,2k}$ from
which the unimodality easily follows. In Section 3, we conjecture
that $B_{n,2k}$ are log-concave for fixed $k$, and shall show that
it holds when $n$ is sufficiently large.

\section{The Unimodality}

In this section, we shall give a proof of the unimodality conjecture of
Palmer, Read and Robinson.
 Let $Q_n$ be the
$n$-dimensional cube represented by a graph whose vertices are
sequences of $1$'s and $-1$'s of length $n$, where two vertices are
adjacent if they differ only at one position. Let $V_n$
denote the set of vertices of $Q_n$, namely,
\[
V_n=\{(\epsilon_1,\epsilon_2,\ldots ,\epsilon_n)\,|\, \epsilon_i=-1\
\mbox{or} \ 1,\ 1\leq i\leq n\}.
\]
By a $2$-coloring of the $Q_n$ we mean an assignment of weights  $1$
or $0$  to the vertices of $Q_n$. The weight
of a $2$-coloring is the sum of weights or the numbers of vertices
with weight 1.
 The center of mass of a coloring $f$ with
$w(f) \neq 0$ is the point whose coordinates are given by
\[
\frac{1}{w(f)}\sum (\epsilon_1,\epsilon_2,\ldots ,\epsilon_n),
\]
where the sum ranges over all black vertices. If $w(f)=0$, we take the
center of mass to be the origin. A $2$-coloring is balanced if its
center of mass coincides with the origin. A pair of vertices of the $n$-cube is
called an antipodal pair if it is of the form $(v, -v)$.  A
$2$-coloring is said to be antipodal if any vertex $v$ and its
antipodal have the same color.

The key idea of our proof  relies on the
following further classification of the set  $\mathcal{B}_{n,2k}$ of
balanced $2$-colorings.

\begin{theorem}\label{le}
Let $\mathcal {B}_{n,2k,i}$ denote the set of the balanced $2$-colorings
in $\mathcal{B}_{n,2k}$ containing exactly $i$ antipodal pairs of
black vertices. Then we have
\begin{equation}\label{main1}
(2^{n-1}-2k+i)|\mathcal{B}_{n,2k,i}|=(i+1)|\mathcal{B}_{n,2k+2,i+1}|,
\end{equation}
for $0\leq i\leq k$ and $1\leq k\leq 2^{n-2}-1$.
\end{theorem}

\proof  We aim to show
that both sides of (\ref{main1}) count the number of  ordered pairs
$(F,G)$, where $F\in \mathcal {B}_{n,2k,i}$ and $G\in \mathcal
{B}_{n,2k+2,i+1}$, such that $G$ can be obtained by changing a pair of
antipodal white vertices of $F$ to black vertices. Equivalently, $F$
can be obtained from $G$ by changing a pair of antipodal black
vertices to white vertices.

First, for each $F\in \mathcal {B}_{n,2k,i}$, we wish to obtain $G$
in $\mathcal {B}_{n,2k+2,i+1}$ by changing a pair of antipodal
 white vertices to black. By the definition of
$\mathcal{B}_{n,2k,i}$, for each $F$ there are $i$ antipodal pairs
of black vertices and $2k-2i$ black vertices whose antipodal
vertices are colored by white. Since  $k\leq 2^{n-2}-1$, that is,
$2^{n-1}-2(k-i)-i>0$, there are exactly  $2^{n-1}-2(k-i)-i$ antipodal pairs
of white vertices in $F$. Thus from each $F\in \mathcal
{B}_{n,2k,i}$, we can obtain $2^{n-2}-2k+i$ different $2$-coloring in
$\mathcal {B}_{n,2k+2,i+1}$ by changing a pair of antipodal white
vertices of $F$ to black. Hence the number of ordered pair $(F,G)$
equals $(2^{n-1}-2k+i)|\mathcal{B}_{n,2k,i}|$.

On the other hand, for each $G\in \mathcal{B}_{n,2k+2,i+1}$, since
there are $i+1$ antipodal pairs of black vertices in $G$, we see
that from $G$ we can obtain $i+1$ different $2$-colorings in
$\mathcal{B}_{n,2k,i}$ by changing a pair of antipodal black
vertices to white. So the number of ordered pairs $(F,G)$ equals
$(i+1)|\mathcal{B}_{n,2k+2,i+1}|$. This completes the proof.\qed

We are ready to prove the unimodality conjecture.

\begin{theorem}\label{main}
For $n\geq 1$, the sequence $\{B_{n,2k}\}_{0\leq k\leq 2^{n-1}}$ is
strictly unimodal with the maximum attained at $k=2^{n-1}$.
\end{theorem}

\proof It is easily seen that  $\{B_{n,2k}\}_{1\leq k\leq
2^{n-1}}$ is symmetric  for any $n\geq 1$. Given a
balanced coloring of the $n$-cube, if we exchange the colors on all
 vertices,  the complementary coloring  is still
balanced. Thus it is sufficient to prove $B_{n,2k}<B_{n,2k+2}$ for
$0\leq k\leq 2^{n-2}-1$.

Clearly, for each $F\in \mathcal{B}_{n,2k}$, there are at most $k$
antipodal pairs of black vertices. It follows that
\[
B_{n,2k}=\sum \limits_{i=0}^{k} |\mathcal {B}_{n,2k,i}|.
\]
We wish to establish the inequality
\begin{equation}\label{le2}
|\mathcal{B}_{n,2k,i}|<|\mathcal{B}_{n,2k+2,i+1}|.
\end{equation}
If it is true,  then

\[
{B}_{n,2k}=\sum \limits_{i=0}^{k} |\mathcal {B}_{n,2k,i}| <\sum
\limits_{i=1}^{k+1} |\mathcal {B}_{n,2k+2,i}| <\sum
\limits_{i=0}^{k+1} |\mathcal {B}_{n,2k+2,i}|={B}_{n,2k+2},
\]
for $0\leq k\leq 2^{n-2}-1$, as claimed in the theorem. Thus it
remains to prove (\ref{le2}). Since $1\leq k\leq 2^{n-2}-1$, it is
clear that
\[
(2^{n-1}-2k+i)-(i+1)=2^{n-1}-2k-1\geq 1.
\]
Applying Theorem \ref{le}, we find that
\[
|\mathcal{B}_{n,2k,i}|< |\mathcal{B}_{n,2k+2,i+1}|,
\]
for $0\leq i\leq k$ and $1\leq k\leq 2^{n-2}-1$, and hence
(\ref{le2}) holds. This completes the proof. \qed

\section{The log-concavity  for fixed $k$}

Log-concave sequences and polynomials often arise in combinatorics,
algebra and geometry, see£» for example, Brenti \cite{brenti1989} and Stanley \cite{stanley1989}.
While $\{B_{n,2k}\}_{k}$  is not log-concave in general,  we shall
show that
 it is log-concave for fixed $k$ and sufficiently large $n$, and we conjecture
 that the log-concavity holds for any given $k$.

\begin{conjecture}\label{con1}
When $0\leq k\leq 2^{n-1}$, we have
\[
B_{n,2k}^2\geq B_{n-1,2k}B_{n+1,2k}.
\]
\end{conjecture}

Palmer, Read and Robinson \cite{palmer} have shown that
\[B_{n,2}=2^{n-1}\] and
\[B_{n,4}=\frac{1}{4^n}((4!)^{n-1}-2^{3n-3}).\]
It is easy to verify  that the sequences $\{B_{n,2}\}_{n\geq
1}$ and $\{B_{n,4}\}_{n\geq 2}$ are both log-concave. Thus in the
remaining of this paper, we shall be concerned only with the case $k\geq 3$.
To be more specific, we shall show that  Conjecture \ref{con1} is true when $n$ is
sufficiently large. Our proof utilizes  the  well-known Bonferroni
inequality,  which can be stated as follows.  Let $P(E_i)$ be the
probability of the event $E_i$, and let  $P\left(\bigcup \limits_{i=1}^n
E_i\right)$ be the probability that at least one of the events $E_1, E_2,\ldots
, E_n$ will occur. Then
\[
P\left(\bigcup _{i=1}^n E_i\right)\leq \sum \limits_{i=1}^{n}P(E_i).
\]

Before we present the proof of the asymptotic log-concavity of the
sequence $\{B_{n, 2k}\}$ for fixed $k$, let us introduce the
$(0,1)$-matrices associated with a balanced $2$-coloring of the
$n$-cube with $2k$ vertices having weight $1$. Since such a
$2$-coloring is uniquely determined by the set of vertices having
weight $1$, we may represent a $2$-coloring by these vertices with
weight $1$. This leads us to consider the set $\mathcal{M}_{n,2k}$
of  $n\times 2k$ matrices such that each row contains $k$ $+1$'s and
$k$ $-1$'s without two identical columns.
 Let $M_{n,2k}=|\mathcal{M}_{n,2k}|$.
It is clear that
 \[ M_{n,2k}=(2k)!B_{n,2k}.\]
Hence the log-concavity of the sequence $\{M_{n,2k}\}_{n\geq \log_2
k+1}$ is equivalent to the log-concavity of the sequence
$\{B_{n,2k}\}_{n\geq \log_2 k+1}$.

Canfield, Gao, Greenhill, McKay and Ronbinson \cite{can} obtained the following
estimate.

\begin{theorem}\label{can1}
If $0\leq k\leq o(2^{n/2})$, then
\[
M_{n,2k}={{2k}\choose
k}^n\left(1-O\left(\frac{k^2}{2^n}\right)\right).
\]
\end{theorem}


To prove the asymptotic log-concavity of $M_{n,2k}$ for fixed $k$,
we need the following  result  that is a stronger property than
Theorem \ref{can1}.


\begin{theorem}\label{dec}
Let $c_{n,k}$ be the real number such that
\begin{equation}\label{m4}
M_{n,2k}={{2k}\choose
k}^n\left(1-c_{n,k}\left(\frac{k^2}{2^n}\right)\right).
\end{equation}
Then we have
\[
c_{n,k}>c_{n+1,k},
\]
when $k\geq 3$ and $n$ is sufficiently large.
\end{theorem}

\proof Let $\mathcal{L}_{n,2k}$ be the set of matrices with every
row consisting of $k$ $-1$'s and $k$ $+1$'s that do not belong to
$\mathcal{M}_{n,2k}$ and $L_{n,2k}=|\mathcal{L}_{n,2k}|$. In other words,
any matrix in $\mathcal{L}_{n,2k}$ has two identical columns. Since the
number of  $n\times 2k$ matrices with each row consisting of $k$
$+1$'s and $k$ $-1$'s eqauls ${{2k}\choose k}^n$. From
(\ref{m4}) it is easily checked  that
\begin{equation}\label{m1}
L_{n,2k}=c_{n,k}\frac{k^2}{2^n}{{2k}\choose k}^n.
\end{equation}

We now proceed to give an upper bound on the cardinality of
$\mathcal{L}_{n+1,2k}$. For each $M\in
\mathcal{L}_{n+1,2k}$, it is easy to see that
the matrix $M'$
obtained from $M$ by deleting the $(n+1)$th row contains
two identical columns as well. Therefore,
every matrix in $\mathcal{L}_{n+1,2k}$ can be obtained from a matrix
in $\mathcal{L}_{n,2k}$ by adding a suitable
 row to a matrix in $\mathcal{L}_{n,2k}$ as the $(n+1)$-th row.
This observation enables us to construct three classes of matrices $M$ from
$\mathcal{L}_{n+1,2k}$ by the properties of $M'$. It is obvious that any matrix in $\mathcal{L}_{n+1, 2k}$ belongs to one of
these three classes.

{\bf Class 1:} There exist at least three identical columns in $M'$.
For each row of $M'$, the probability that the three prescribed
positions of this row are identical equals
\[
2{{2k-3}\choose k}\left/{{2k}\choose k}\right..
\]
Here the factor $2$ indicates that there are two choices for the
values at the  prescribed positions. Consequently, the probability
that the three prescribed columns in $M'$ are identical equals
\[
\left(2{{2k-3}\choose k}\left/{{2k}\choose
k}\right. \right)^n=\left(\frac{k-2}{2(2k-1)}\right)^n<\frac{1}{4^n}.
\]

By the Bonferroni inequality, the probability that there are at
least three identical columns in $M'$ is bounded by
$\frac{8k^3}{4^n}$. Because the number of  $(n+1)\times 2k$ matrices
with each row consisting of $k$ $+1$'s and $k$ $-1$'s is
${{2k}\choose k}^{n+1}$,
 the number of matrices $M$ in $\mathcal{L}_{n+1,2k}$ with $M'$ containing
 at least three
 identical columns is bounded by
\[
\frac{8k^3}{4^n}{{2k}\choose k}^{n+1}.
\]

{\bf Class 2:} There exist at least two pairs of identical columns
in $M'$. For any two prescribed pairs $(i_1,i_2)$ and $(j_1,j_2)$ of
columns, let us estimate the probability that in $M'$ the $i_1$-th
column is identical to the $i_2$-th column and the $j_1$-th column
is identical to the $j_2$-th column, that is, for any row of $M'$,
the value of the $i_1$-th (respectively, $j_1$-th) position is equal
to the value of the $i_2$-th (respectively, $j_2$-th) position. We
have two cases for each row of $M'$. The first case is that the
values at the positions $i_1$, $i_2$, $j_1$ and $j_2$ are all
identical. The probability for any given row to be in this case
equals
\[
2{{2k-4}\choose {k-4}}\left/{{2k}\choose k}\right. .
\]
Again, the factor $2$ comes from the two choices for the values at the prescribed positions.

The second  case is that the value of the $i_1$-th position is different from the value of the $j_1$-th position. In this case, we have either  the values at the $i_1$-th and $i_2$-th positions are $+1$ and the values at the $j_1$-th and $j_2$-th positions are $-1$ or the values at $i_1$-th and $i_2$-th position are $-1$ and the values at the $j_1$-th and $j_2$-th positions are $+1$. Thus the probability  for
 any given row to be in this case equals
\[
2{{2k-4}\choose
{k-2}}\left/{{2k}\choose k}\right..
\]
Combining the above two case, we see that  when $k\geq 3$, the probability that  $M'$ has
two prescribed pairs of identical columns equals
\[
\left(2{{2k-4}\choose {k-4}}\left/{{2k}\choose k}\right. +2{{2k-4}\choose
{k-2}}\left/{{2k}\choose k} \right. \right)^n<\frac{1}{4^n}.
\]
Again, by the Bonferroni inequality, the probability that there exist
at least two pairs of identical columns of $M'$ is bounded by
$\frac{16k^4}{4^n}$. It follows that the number of matrices $M$ in
$\mathcal{L}_{n+1,2k}$ with $M'$ containing at least two pairs of
identical columns is bounded by
\[
\frac{16k^4}{4^n}{{2k}\choose k}^{n+1}.
\]

{\bf Class 3:} There exists exactly one pair of identical
columns in $M'$. By the definition, the
number of matrices $M'$ containing exactly one pair of identical columns is bounded by
$L_{n,2k}$. On the other hand, it is easy to see that for each $M'$
containing exactly one pair of identical columns, there are
\begin{equation}\label{m6}
2{{2k-2}\choose k}=\frac{k-1}{2k-1}{{2k}\choose k}
\end{equation}
matrices of $\mathcal{L}_{n+1,2k}$ which can be obtained by adding a
suitable row as the $(n+1)$-th row. Combining (\ref{m1}) and
(\ref{m6}), we find that the number of matrices $M$ of $\mathcal{L}_{n+1,2k}$
such  $M'$ contains exactly one pair of
identical columns is bounded by
\[
\frac{k-1}{2k-1}c_{n,k}\frac{k^2}{2^n}{{2k}\choose k}^{n+1}.
\]
Clearly, $L_{n+1,2k}$ is bounded by sum of the cardinalities of the above three classes.
This yields the upper bound
\[
L_{n+1,2k}<\frac{8k^3}{4^n}{{2k}\choose
k}^{n+1}+\frac{16k^4}{4^n}{{2k}\choose
k}^{n+1}+\frac{k-1}{2k-1}c_{n,k}\frac{k^2}{2^n}{{2k}\choose
k}^{n+1},
\]
when $k\geq 3$ and $n$ is sufficiently large.

Now we claim that
\begin{equation}\label{m5} \frac{8k^3}{4^n}+\frac{16k^4}{4^n}=o\left(c_{n,k}\frac{k^2}{2^n}\right).\end{equation} Notice that
the probability that a specified pair of columns in $M'$ are identical is
\[
\left(2{{2k-2}\choose k}\left/{{2k}\choose
k}\right.\right)^n=\left(\frac{k-1}{2k-1}\right)^n.
\]
Since  $c_{n,k}\frac{k^2}{2^n}$ is the
probability  that there exists at least
two identical columns in $M'$, we deduce that
\[ c_{n,k}\frac{k^2}{2^n}>\left(2{{2k-2}\choose k}\left/{{2k}\choose
k}\right. \right)^n=\left(\frac{k-1}{2k-1}\right)^n>\frac{1}{3^n}.\]
But when $n$ is sufficiently large, we have
\[ \frac{8k^3}{4^n}+\frac{16k^4}{4^n}=o\left(\frac{1}{3^n}\right),\]
which implies (\ref{m5}).
Since  $ \frac{k-1}{2k-1}<\frac{1}{2}$,  it follows from (\ref{m5}) that
\begin{equation}\label{m2}
L_{n+1,2k}<c_{n,k}\frac{k^2}{2^{n+1}}{{2k}\choose k}^{n+1},
\end{equation}
when $n$ is sufficient large.
Restating formula (\ref{m1}) for $n+1$, we have
\begin{equation}\label{m3}
L_{n+1,2k}=c_{n+1,k}\frac{k^2}{2^{n+1}}{{2k}\choose k}^{n+1}.
\end{equation}
Combining (\ref{m2}) and (\ref{m3}) gives
\[
c_{n,k}>c_{n+1,k},
\]
for sufficiently large $n$. This completes the proof.\qed

Applying Theorem \ref{dec}, we arrive at  the following result.

\begin{theorem}\label{suf}
When $n$ is sufficiently large,
 \[ M_{n,2k}^2>M_{n-1,2k}M_{n+1,2k}.\]
\end{theorem}

\proof We only consider the case $k\geq 3$. Let \[
M_{n,2k}={{2k}\choose k}^n\left(1-c_{n,k}\frac{k^2}{2^n}\right).\]
Then
\begin{eqnarray*}
& &M_{n,2k}^2-M_{n-1,2k}M_{n+1,2k}\\[5pt]
&=&{{2k}\choose
k}^{2n}\left[\left(1-c_{n,k}\frac{k^2}{2^n}\right)^2-
\left(1-c_{n+1,k}\frac{k^2}{2^{n+1}}\right)\left(1-c_{n-1,k}\frac{k^2}{2^{n-1}}\right)\right]\\[5pt]
&=&{{2k}\choose
k}^{2n}\left[-c_{n,k}\frac{k^2}{2^{n-1}}+c_{n,k}^2\frac{k^2}{4^n}+c_{n+1,k}\frac{k^2}{2^{n+1}}
+c_{n-1,k}\frac{k^2}{2^{n-1}}-c_{n-1,k}c_{n+1,k}\frac{k^2}{4^n}\right].
\end{eqnarray*}
By Theorem \ref{dec},  we have $c_{n-1,k}>c_{n,k}$ when $k\geq 3$ and
$n$ is sufficiently large. This implies that
\[ c_{n,k}\frac{k^2}{2^{n-1}}<c_{n+1,k}\frac{k^2}{2^{n-1}},\]
when $k\geq 3$ and $n$ is sufficiently large. By the proof of
Theorem \ref{can1} in \cite{can}, we have $c_{n-1,k}<4$ for any $n$.
It follows that
\[
c_{n-1,k}c_{n+1,k}\frac{k^2}{4^n}<c_{n+1,k}\frac{k^2}{4^{n-1}}\leq
c_{n+1,k}\frac{k^2}{2^{n+1}},
\]
when $n\geq 3$. Hence
\[
M_{n,2k}^2>M_{n-1,2k}M_{n+1,2k},
\]
for sufficiently large $n$.  This completes this proof. \qed

Since $M_{n,2k}=(2k)!B_{n,2k}$, Theorem \ref{suf} implies the
asymptotic log-concavity of $B_{n,2k}$ for fixed $k$.

\begin{corollary}
When $n$ is sufficiently large,  \[B_{n,2k}^2>
B_{n-1,2k}B_{n+1,2k}.\]
\end{corollary}

\vspace{2mm}

\noindent {\bf Acknowledgments.} This work was supported by the 973
Project, the PCSIRT Project of the Ministry of Education,  and the
National Science Foundation of China.

\end{document}